\documentclass{amsart}
\usepackage{amsthm,amscd,amssymb}
\usepackage{diagrams}

\newcommand{\rsp}{\raisebox{0em}[2.7ex][1.3ex]{\rule{0em}{2ex}}}
\newcommand{\F}{\mathbb F}
\newcommand{\N}{\mathbb N}
\newcommand{\Z}{\mathbb Z}
\newcommand{\Q}{\mathbb Q}
\newcommand{\OO}{\mathcal O}
\newcommand{\cC}{\mathcal C}

\newcommand{\Cl}{\operatorname{Cl}}

\newcommand{\disc}{\operatorname{disc}}
\newcommand{\Am}{\operatorname{Am}}
\newcommand{\Aut}{\operatorname{Aut}}
\newcommand{\rk}{\operatorname{rk}\,}
\newcommand{\st}{{\operatorname{st}}}
\newcommand{\Gal}{\operatorname{Gal}}
\newcommand{\Spl}{\operatorname{Spl}}
\newcommand{\la}{\langle}
\newcommand{\ra}{\rangle}
\newcommand{\too}{\longmapsto}
\newcommand{\lra}{\longrightarrow}
\newcommand{\hra}{\hookrightarrow}
\newcommand{\fa}{\mathfrak a}
\newcommand{\fb}{\mathfrak b}

\newcommand{\eps}{\varepsilon}

\renewcommand{\theta}{\vartheta}

\newtheorem{thm}{Theorem}
\newtheorem{prop}{Proposition}
\newtheorem{lem}{Lemma}

\title{Class Groups of Dihedral Extensions}
\author{Franz Lemmermeyer}
\address{CSU San Marcos, Dept. Mathematics \\
333 S Twin Oaks Valley Rd \\
San Marcos, CA 92096-0001 \\ USA}
\email{franzl@csusm.edu}

\begin{document}

\begin{abstract}
Let $L/F$ be a dihedral extension of degree $2p$, where $p$ is
an odd prime. Let $K/F$ and $k/F$ be subextensions of $L/F$ with 
degrees $p$ and $2$, respectively. Then we will study relations
between the $p$-ranks of the class groups $\Cl(K)$ and $\Cl(k)$.
\end{abstract}

\maketitle

\section{A Short History of Reflection Theorems}

Results comparing the $p$-rank of class groups of different
number fields (often based on the interplay between Kummer
theory and class field theory) are traditionally called 
`reflection theorems'; the oldest such result is due to 
Kummer himself: let $h^+$ and $h^-$ denote the plus and
the minus $p$-class number of $K = \Q(\zeta_p)$, respectively;
then Kummer observed that $p \mid h^+$ implies $p \mid h^-$,
and this was an important step in verifying Fermat's Last Theorem
(that is, checking the regularity of $p$) for exponents $< 100$.
Kummer's result was improved by Hecke \cite{Hecke} (see also Takagi 
\cite{Tak}):

\begin{prop}
Let $p$ be an odd prime, $k = \Q(\zeta_p)$, and let
$\Cl_p^+(k)$ and $\Cl_p^-(k)$ denote the plus and the
minus part of $\Cl_p(k)$. Then $\rk \Cl_p^+(k) \le 
\rk \Cl_p^-(k)$.
\end{prop}

Analogous inequalities hold for the eigenspaces of the 
class group $\Cl(k)$ under the action of the Galois group;
see e.g. \cite{IR}.

Scholz \cite{Sch32} and Reichardt \cite{Rei34} discovered a 
similar connection between the $3$-ranks of class groups of 
certain quadratic number fields:

\begin{prop}\label{PS}
Let $k^+ = \Q(\sqrt{m}\,)$ with $m \in \N$, and put 
$k^- = \Q(\sqrt{-3m}\,)$; then the $3$-ranks $r_3^+$ and 
$r_3^-$ of $\Cl(k^+)$ and $\Cl(k^-)$ satisfy the inequalities 
$r_3^+ \le r_3^- \le r_3^+ + 1$.
\end{prop}

Leopoldt \cite{Leo58} later generalized these propositions 
considerably and called his result the ``Spie\-ge\-lungs\-satz''. 
For expositions and generalizations, see Kuroda \cite{Kur}, 
Oriat \cite{Ori76,Ori79}, Satg\'e \cite{Sat76}, Oriat \& Satg\'e 
\cite{OS79}, and G.~Gras \cite{GrasTR}.

Damey \& Payan \cite{DP70} found an analog of Proposition
\ref{PS} for $4$-ranks of class groups of quadratic number fields:

\begin{prop}
Let $k^+ = \Q(\sqrt{m}\,)$ be a real quadratic number field, 
and put $k^- = \Q(\sqrt{-m}\,)$. Then the $4$-ranks $r_4^+$ 
and $r_4^-$ of $\Cl^+(k^+)$ (ideal class group in the strict 
sense) and $\Cl(k^-)$ satisfy the inequalities
$r_4^+ \le r_4^- \le r_4^+ + 1$.
\end{prop}

Other proofs were given by G.~Gras \cite{Gras3}, Halter-Koch \cite{HK84}, 
and Uehara \cite{Ueh89}; for a generalization, see Oriat \cite{Ori77a,Ori77b}. 

In 1974, Callahan \cite{Cal} discovered the following result; 
although it gives a connection between $p$-ranks of class groups
of different number fields, its proof differs considerably
from those of classical reflection theorems:

\begin{prop}
Let $k$ be a quadratic number field with discriminant $d$,
and suppose that its class number is divisible by $3$. Let
$K$ be one of the cubic extensions of $\Q$ with discriminant
$d$ (then $Kk/k$ is a cyclic unramified extension of $k$), 
and let $r_3(k)$ and $r_3(K)$ denote the $3$-ranks of $\Cl(k)$
and $\Cl(K)$, respectively. Then $r_3 = r_2 -1$.
\end{prop}

Callahan could only prove that $r_2-2 \le r_3 \le r_2 -1$, 
but conjectured that in fact $r_3 = r_2 -1$. This was verified
later by  G.~Gras \cite{Gras} and Gerth \cite{Ger}. Callahan's
result was generalized by B\"olling \cite{Boe}: 

\begin{prop}\label{PCal}
Let $L/\Q$ be a normal extension with Galois group the dihedral 
group of order $2p$, where $p$ is an odd prime, and let $K$
be any of its subfields of degree $p$. Assume that the quadratic
subfield $k$ of $L$ is complex. Then 
$$ r_p(k) - 1 \le r_p(K) \le \frac{p-1}2 (r_p(k)-1), $$
where $r_p(k)$ and $r_p(K)$ denote the $p$-ranks of the class groups
of $k$ and $K$, respectively.
\end{prop}

It is this result that we generalize to arbitrary base fields 
in this article. Our proof will be much less technical than
B\"olling's, who used the Galois cohomological machinery presented 
in Koch's book \cite{Koch}.

We conclude our survey of reflection theorems with the following
result by Kobayashi \cite{Kob73} (see also Gerth \cite{Ger75}):

\begin{prop}\label{PKob}
Let $m$ be a cubefree integer not divisible by any prime 
$p \equiv 1 \bmod 3$, and put $K = \Q(\sqrt[3\,]{m}\,)$ and 
$L = K(\sqrt{-3}\,)$. Then $\rk \Cl_3(L) = 2 \cdot \rk \Cl_3(K)$. 
\end{prop}

This was generalized subsequently by G.~Gras \cite{Gras} to
the following result; $\Spl(k/\Q)$ denotes the set of primes 
in $\Q$ that split in $k$.

\begin{prop}
Let $K$ be a cubic number field with normal closure $L$. Assume that
$\Gal(L/\Q) \simeq S_3$, and let $k$ denote the quadratic subfield of
$L$. If no prime $p \in \Spl(k/\Q)$ ramifies in $L/k$, and if 
$3 \nmid h(k)$, then $\rk \Cl_3(L) = 2 \cdot \rk \Cl_3(K)$. 
\end{prop}

It seems plausible that the results of Proposition \ref{PCal} hold 
for a large variety of nonabelian extensions. Computer experiments 
suggest a rather simple result normal extensions of $\Q$ with Galois 
group $A_4$. In fact, consider a cyclic cubic extension $k/\Q$, and 
let $2r$ denote the $2$-rank of $\Cl(k)$. Then there exist $r$ 
nonconjugate quartic extensions $K/\Q$ such that (cf. \cite{H})
\begin{enumerate}
\item $Kk$ is the normal closure of $K/\Q$, and
        $\Gal(Kk/\Q) \simeq A_4$;
\item $Kk/k$ is an unramified normal extension with 
        $\Gal(Kk/k) \simeq (2,2)$. 
\end{enumerate}

\medskip\noindent
{\bf Conjecture 1.}
{\em Let $L/\Q$ be an $A_4$ extension unramified over its cubic 
subfield $k$, and let $K$ denote one of the four conjugate quartic 
subfields of $L$. Then we have the following inequalities:}
\begin{eqnarray*}
     \rk \Cl_2(k) & \ge \ \rk \Cl_2(K)   & \ge \ \rk \Cl_2(k)-2, \\
 \rk \Cl_2^+(k)+1 & \ge \ \rk \Cl_2^+(K) & \ge \ \rk \Cl_2^+(k)-1
\end{eqnarray*}
Examples show that these inequalities are best possible. 
In fact, consider the cyclic cubic extension $k$ generated 
by a root of the cubic polynomial $f(x) = x^3-ax^2-(a+3)x-1$; 
choose $b \in \N$ odd and $a = \frac12(b^2-3)$. Then 
$L = k(\sqrt{\alpha-2},\sqrt{\alpha'-2}\,)$ is an $A_4$-extension
of $\Q$, $L/k$ is unramified, and using PARI we find that
$\Cl(k) \simeq (4,4,2,2)$, $\Cl(K) \simeq (2,2)$ for $a = 143$, and
$\Cl(k) \simeq (114,2)$, $\Cl(K) \simeq (4,2)$ for $a = 1011$. 

\section{The Main Results}

Let $p = 2m+1$ be an odd prime and let
$$ D_p = \la \sigma, \tau: \sigma^p = \tau^2 = 1, 
             \tau\sigma\tau = \sigma^{-1} \ra$$
denote the dihedral group of order $2p$. Put 
$\nu = 1 + \sigma + \sigma^2 + \ldots + \sigma^{p-1}$;
then a simple calculation gives $\nu \tau = \tau \nu$
and $\nu \sigma = \sigma \nu$.

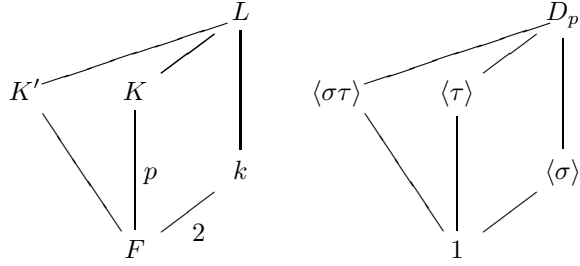
$$\hbox{
\begin{diagram}[width=2em,height=1.5em]
& &      &      &  L \\
&  &  &  \ruLine(4,2) \ruLine  \ruLine &  \dLine  \\
K'\ & &  K   &      &    \\
& \rdLine(2,4)&      &      &    \\
& &  \dLine_p  &      &  k \\
& &      & \ruLine_2 &    \\
& &   F  &      &  
\end{diagram} \qquad 
\begin{diagram}[width=2em,height=1.5em]
& &      &      &  D_p \\
&  &  & \ruLine(4,2) \ruLine &  \dLine  \\
\la \sigma\tau \ra\phantom{lP} & &  \la \tau \ra   &      &    \\
& \rdLine(2,4)&      &      &    \\
& &  \dLine  &      &  \la \sigma \ra \\
& &      & \ruLine &    \\
& &   1  &      &  
\end{diagram}}
$$

Let $F$ be a number field with class number not divisible by $p$,
$L/F$ a dihedral extension with Galois group $\Gal(L/F) \simeq D_p$,
$k$ its quadratic subfield, and $K$ the fixed field of $\tau$.
Note that $K' = K^{\sigma^m}$ is the fixed field of $\sigma\tau$.

In the sequel, $L/k$ will always be unramified. Our main result
generalizes B\"olling's theorem to base fields $F$ with class
number prime to $p$:

\begin{thm}\label{T1}
Let $F$ be a number field with class number not divisible by $p$,
let $L/F$ be a $D_p$-extension such that $L/k$ is unramified, and 
let $r_p(k)$ and $r_p(K)$ denote the $p$-ranks of the class groups
of $k$ and $K$, respectively; then
$$ r_p(k) - 1 - e \le r_p(K), $$
where $p^e = (E_F:N E_K)$. 
\end{thm}

The idea of the proof of Theorem \ref{T1} is to compare the 
class groups of $K$ and $k$ by lifting them to $L$ and studying 
homomorphisms between certain subgroups of $\Cl_p(L)$. 
We get inequalities for the ranks by computing the orders 
of elementary abelian $p$-groups.

If $F = \Q$ or if $F$ is a complex quadratic number field (different
from $\Q(\sqrt{-3}\,)$ if $p = 3$), then $e = 0$ since in these
cases $E_F$ is torsion of order not divisible by $p$. 

Actually, B\"olling's upper bound from Prop. \ref{PCal} is 
conjectured to be valid in general:

\medskip\noindent
{\bf Conjecture 2.} 
{\em Under the assumptions of Theorem \ref{T1} we have}
$$ r_p(K) \le \frac{p-1}2 (r_p(k)-1). $$
\medskip

We will prove Conjecture 2 if $p = 3$, if $r_p(k) = 1$, or if 
$\Cl_p(k) = (p,p)$; by B\"olling's result, the upper bound holds 
if $F = \Q$ and $k$ is complex quadratic.

\medskip\noindent
{\bf Conjecture 3.} 
{\em Fix an odd prime $p$ and a number field $F$ with
  class number not divisible by $p$. Then for every
integer $e$ with $0 \le e \le \dim E_F/E_F^p$,
every integer $r \ge 1$  and every $R \ge 0$ such that
$r-1-e \le R \le \frac{p-1}2 (r-1)$ there exist dihedral 
extensions $L/F$ satisfying the assumptions of Theorem \ref{T1}
such that $r_p(k) = r$, $r_p(K) = R$, and $(E_F:N_{K/F}E_K) = p^e$.}
\medskip

A proof of Conjecture 3 seems to be completely out of reach; it 
expresses the expectation that the bounds in Theorem \ref{T1} and 
Conjecture 2 are best possible.

Using the results needed for the proof of Theorem \ref{T1}, 
we get the following class number formula almost for free:

\begin{thm}\label{T2}
Let $L/F$ be a dihedral extension of degree $2p$, where $p$ is
an odd prime, and assume that $L$ is unramified over the quadratic
subextension $k$ of $L/K$. Let $q = (E_L:E_K E_{K'} E_k)$ denote 
the unit index of $L/F$ and write $a = 1 + \lambda(k) - \lambda(F)$, 
where $\lambda(M)$ denotes the $\Z$-rank of the unit group of a 
number field $M$. Then 
\begin{equation}\label{ECNF}
   h_L = p^{-a} q h_k \Big(\frac{h_K}{h_F}\Big)^2. 
\end{equation}
\end{thm}

In the special case $F = \Q$, an arithmetic proof of the
class number formula for dihedral extensions of degree $2p$
(even without the restriction that $L/k$ be unramified)
was given by Halter-Koch \cite{HK}.

A simple application of the lower bound in Theorem \ref{T1}
gives 

\begin{thm}\label{TNo}
Let $L/F$ be as in Theorem \ref{T1}. If $r_p(k) \ge e+2$, then there 
exists a normal unramified extension $M/k$ (containing $L$) with 
$\Gal(M/k) \simeq E(p^3)$, the nonabelian group of order $p^3$ and 
exponent $p$.
\end{thm}

In the special case where $F$ is $\Q$ or a complex quadratic number
field $\ne \Q(\sqrt{-3}\,)$ this was proved by Nomura \cite{Nom}
(note that $e = 0$ in these cases).

\section{Preliminaries}
In this section we collect some results that will be needed 
in the sequel. 

Let $\Am = \Am(L/k) = \{c \in \Cl_p(L) : c^\sigma = c\}$ denote 
the ambiguous $p$-class group and 
$\Am_\st = \{c = [\fa] \in \Am: \fa^\sigma = \fa\}$ its subgroup
of strongly ambiguous ideal classes. Since $L/k$ is unramified,
ambiguous ideals are ideals from $k$, hence $\Am_\st = \Cl_p(k)^j$, 
where $j: \Cl(k) \lra \Cl(L)$ is the transfer of ideal classes. 
This proves

\begin{lem}
For unramified extensions $L/k$, the sequence
$$ \begin{CD}\label{Ekap}
   1 @>>> \kappa_{L/k} @>>> \Cl_p(k) @>j>> \Am_\st @>>> 1,
    \end{CD} $$
where $\kappa_{L/k}$ is the capitulation kernel, is exact.
\end{lem}

The next lemma is classical; it measures the difference between
the orders of ambiguous and strongly ambiguous ideal classes:

\begin{lem}
Let $L/k$ be a cyclic extension of prime degree $p$. Then the 
factor group $\Am/\Am_\st$ is an elementary abelian $p$-group.
In fact, we have the exact sequence
\begin{equation}\label{E2}
 \begin{CD}
   1 @>>> \Am_\st @>>> \Am @>{\theta}>> E_k \cap NL^\times/NE_L @>>> 1. 
 \end{CD} \end{equation}
If, in addition, $L/k$ is unramified, then $E_k \cap NL^\times = E_k$, 
and we find
\begin{equation}\label{ECam}
 \begin{CD} 1 @>>> \Am_\st @>>> \Am @>{\theta}>> E_k/NE_L @>>> 1. \end{CD} 
\end{equation}
\end{lem}

\begin{proof}
Let us start by defining $\theta$. Write $c = [\fa] \in \Am$; then 
$\fa^{\sigma-1} = (\alpha)$ and $\eps := N\alpha \in E_k$. Now put
$\theta(c) = \eps NE_L$. We claim that $\theta$ is well defined:
in fact, if $c = [\fb]$ and $\fb^{\sigma-1} = (\beta)$, then 
$\fa = \gamma \fb$ for some $\gamma \in L^\times$; thus
$\fa^{\sigma-1} = \gamma^{\sigma-1} \fb^{\sigma-1}$, and this
shows that $\alpha = \eta \gamma^{\sigma-1} \beta$, hence
$N\alpha = N\eta N\beta$, and therefore 
$N\alpha \cdot NE_L = N\beta \cdot NE_L$.

We have $c \in \ker \theta$ if and only if $\eps = N_{L/k} \eta$
for some unit $\eta$. Then $N(\alpha/\eta) = 1$, hence
$\alpha/\eta = \beta^{1-\sigma}$, therefore $\beta\fa$ is an 
ambiguous ideal, and this implies that $c \in \Am_\st$.
Conversely, if $c \in \Am_\st$, then $c = [\fa]$ with
$\fa^{\sigma-1} = (1)$, hence $\theta(c) = 1$.

It remains to show that $\theta$ is surjective. Given
$\eps \in E_k \cap NL^\times$, write $\eps = N_{L/k}\alpha$ 
for some $\alpha \in L^\times$: then $N_{L/k}(\alpha) = (1)$, 
hence Hilbert's theorem 90 for ideals implies that 
$(\alpha) = \fa^{\sigma-1}$ for some ideal $\fa$ in $\OO_L$; 
clearly $\eps N E_L  = \theta([\fa])$, and this proves the claim.

Finally we have to explain why $E_k \cap NL^\times = E_k$ if
$L/k$ is unramified. In this case, every unit is a local norm 
everywhere (in the absence of global ramification, every local 
extension is unramified, and units are always norms in unramified 
extensions of local fields), hence a global norm by Hasse's 
norm residue theorem for cyclic extensions $L/k$. 
\end{proof}

\begin{prop}[Furtw\"angler's Theorem 90]
If $L/k$ is a cyclic unramified extension of prime degree $p$ and 
$\Gal(L/k) = \la \sigma \ra$, then $\Cl_p(L)[N] = \Cl_p(L)^{1-\sigma}$.
\end{prop}

\begin{proof}
This is a special case of the principal genus theorem of 
classical class field theory; see \cite{PGT}.
\end{proof}

\begin{lem}\label{LST}
Let $A$ be a $D_p$-module; then $A^{1-\sigma} = A^{1+\tau}A^{1+\sigma\tau}$.
\end{lem}

\begin{proof}
For $a \in A$ we have 
$a^{1-\sigma} = (a^{-\sigma})^{1+\tau} a^{1+\sigma\tau}$.
\end{proof}

\section{Galois Action}

Let $m > 1$ be an integer, $p \equiv 1 \bmod m$ an odd prime,
and $r$ an element of order $m$ in $(\Z/p\Z)^\times$. 
Consider the Frobenius group
$$ F_{mp} = \la \sigma, \tau: \sigma^p = \tau^m = 1,
                 \tau^{-1}\sigma \tau = \sigma^r \ra. $$ 
In the following, let $s \in \F_p$ denote the inverse of $r$;
note that $\tau \sigma \tau^{-1} = \tau^s$.

Let $A$ be an abelian $p$-group and $F_{mp}$-module. 
Then the action of $H = \la \tau \ra \simeq \Z/m\Z$ allows
us to decompose $A$ into eigenspaces 
$$ A = \bigoplus_{j = 0}^{m-1} A(j), $$
where $A_j = \{a \in A: a^\tau = a^{r^j}\}$. Note that
$A(j) = e_jA$ for
$$ e_j = \frac1m \sum_{k=0}^{m-1} s^{jk} \tau^k; $$ 
the set $\{e_0, e_1, \ldots, e_{m-1} \}$ is a complete set of
orthogonal idempotents of the group ring $(\Z/m\Z)[F_{mp}]$. 
Also observe that $A(0) = A^H$ is the fixed module of $A$ 
under the action of $H$.

\begin{lem}
Let $A, B, C$ be abelian $p$-groups and $H$-modules. If
$$ \begin{CD} 1 @>>> A @>\iota>> B @>\pi>> C @>>> 1 \end{CD} $$
is an exact sequence of $H$-modules, then so is
$$ \begin{CD} 1 @>>> A(j) @>>> B(j) @>>> C(j) @>>> 1 \end{CD} $$
for every $0 \le j \le m-1$.
\end{lem}

\begin{proof}
This is a purely formal verification based on the fact that,
by the assumption that $\iota$ and $\pi$ be $H$-homomorphisms, 
the homomorphisms $e_j$ commute with $\iota$ and $\pi$.
\end{proof}

\begin{prop}\label{PES}
Assume that $A$, $B$, $C$ are abelian $p$-groups and $H$-modules, 
that 
$$ \begin{CD} 1 @>>> A @>{\iota}>> B @>{\pi}>> C @>>> 1 \end{CD} $$
is an exact sequence of abelian groups, and that 
$\iota(a^\tau) = \iota(a)^\tau$ and $\pi(b^\tau) = \pi(b)^{s\tau}$. 
Then 
$$ \begin{CD} 1 @>>> A(j) @>{\iota}>> B(j) @>{\pi}>> C(j+1) @>>> 1 
  \end{CD} $$  
is exact for every $0 \le j \le m-1$.
\end{prop}

\begin{proof}
Define an $H$-module $C'$ by putting $C = C'$ as an abelian group
and letting $\tau$ act on $C'$ via $c \too c^{s\tau}$. Then
$$ \begin{CD} 1 @>>> A @>{\iota}>> B @>{\pi}>> C' @>>> 1 \end{CD} $$
is an exact sequence of $H$-modules, and taking the $e_j$-part we
get the exact sequence
$$ \begin{CD} 1 @>>> A(j) @>>> B(j) @>>> C'(j) @>>> 1. \end{CD} $$
But $C'(j) = \{c \in C: c^{s\tau} = c^{r^j}\}
           = \{c \in C: c^\tau = c^{r^{j+1}}\} = C(j+1)$.
\end{proof}

Before we can apply the results above to our situation, we have
to check that the homomorphism $\theta$ in (\ref{E2}) satisfies
the assumption of Prop. \ref{PES}.

\begin{lem} 
Let $L/F$ be an $F_{mp}$-extension. Then the map $\theta$ in 
$(\ref{E2})$ (and therefore also in $(\ref{ECam})$) has the property 
$\theta(c^\tau) = \theta(c)^{s\tau}$.
\end{lem}

\begin{proof}
Write $c = [\fa]$, $\fa^{\sigma-1} = (\alpha)$, and 
$N_{L/k} \alpha = \eps$; then $\theta(c) = \eps N_{L/k} E_L$. 
We have $\tau(\sigma-1) = (\sigma^s-1)\tau = (\sigma-1)\phi \tau$ for 
$\phi = 1 + \sigma + \ldots + \sigma^{s-1}$, hence
$(\fa^\tau)^{\sigma-1} = (\fa^{\sigma-1})^{\phi \tau}
  = (\alpha^{\phi\tau})$. Since the norm 
$1 + \sigma + \ldots + \sigma^{p-1}$ is in the center
of $\Z[F_{mp}]$, we get
$N_{L/k} (\alpha^{\phi\tau}) = (N_{L/k} \alpha)^{\phi\tau}
 = \eps^{s\tau}$, and this shows $\theta(c^\tau) = c^{s\tau}$ 
as claimed.
\end{proof}

Let us now specialize to the case $m=2$, where $F_{2p} = D_p$ 
is the dihedral group of order $2p$. For $D_p$-modules $A$ we put
$$ A^+ = A(0) = \{a \in A: a^\tau = a\} \quad \text{and} \quad 
   A^- = A(1) = \{a \in A: a^\tau = a^{-1}\}. $$
If $A$ is finite and has order coprime to $p$, then 
$A = A^+ \oplus A^-$, $A^+ = A^{1+\tau}$ and $A^- = A^{1-\tau}$. 
In the following, let $H = \la \tau \ra$ denote the subgroup of 
$D_p$ generated by $\tau$.

The main ingredient in our proof of Theorem \ref{T1} will be
the following result, which was partially proved by G.~Gras \cite{Gras}:

\begin{thm}
Let $L/F$ be as above, and assume in particular that $L/k$ is
unramified. Then there is an exact sequence 
\begin{equation}\label{EES}
\begin{CD}
    1 @>>> \Am_\st @>{\iota}>>  \Am^- @>\theta>> E_F/NE_K @>>> 1.
   \end{CD} \end{equation}
Moreover, $\Am^+ \simeq (E_k/NE_L)^-$; in particular,
$\Am^+$ is an elementary abelian group of order 
$p^{\rho-1-e}$, where $p^\rho = \# \kappa_{L/k}$ is the 
order of the capitulation kernel and $p^e = (E_F:NE_K)$.
\end{thm}

\begin{proof}
We apply Proposition \ref{PES} with $i = 1$ to (\ref{ECam}). 
Clearly $\tau$ acts as $-1$ on $\Am_\st$, hence $\Am_\st^- = \Am_\st$.
Thus we only have to show that the plus part of $E_k/N_{L/k}E_L$
is isomorphic to $E_F/N_{K/F} E_K$. By sending $\eps N_{K/F} E_K$
to $\eps N_{L/k}E_L$ we get a homomorphism $\psi:E_F/N_{K/F} E_K \lra 
(E_k/N_{L/k}E_L)^+$. 

We claim that $\psi$ is injective; in fact, 
$\ker \psi = \{\eps N_{K/F} E_K: \eps \in  N_{L/k}E_L\}$; but 
$\eps = N_{L/k} \eta$ implies $\eps^2 = \eps^{1+\tau} = 
   N_{L/k} \eta^{1+\tau} = N_{K/F} \eta^{1+\tau}$.
Thus $\eps^2 \in N_{K/F} E_K$, hence so is 
$\eps^{1+p} = (\eps^2)^{(p+1)/2}$. Since 
$E_F/N_{K/F} E_K$ is a $p$-group, we have 
$\eps \in N_{K/F} E_K$.

Moreover, $\psi$ is surjective: in fact, if 
$\eps N_{L/k}E_L$ is fixed by $\tau$, then 
$\eps^2 N_{L/k}E_L = \eps^{1+\tau} N_{L/k}E_L$ is clearly
in the image of $\psi$, and the claim follows again from the
fact that $E_k/N_{L/k} E_L$ is a $p$-group.

Applying Proposition \ref{PES} with $i = 0$ yields 
the isomorphism  $\Am^+ \simeq (E_k/NE_L)^-$; since
$E_k/NE_L$ is elementary abelian, so is $\Am^+$. 
Moreover, the decomposition into eigenspaces 
$E_k/NE_L = (E_k/NE_L)^- \oplus (E_k/NE_L)^+$ shows
$$ \# \Am^+ = \frac{(E_k:NE_L)}{(E_F:NE_K)}.$$

The exact sequence (\ref{Ekap})
shows that $p^\rho = \# \Cl_p(k)/\# \Am_\st$; since 
$$ \# \Am_\st = \frac{\# \Cl_p(k)}{p(E_k:NE_L)}, $$
this implies that $p^\rho = p(E_k:NE_L)$, hence
$\# \Am^{1+\tau} =  \frac{(E_k:NE_L)}{p(E_F:NE_K)} = p^{\rho-1-e}$.
\end{proof}

\section{The Class Number Formula}

As a simple application of the exact sequence (\ref{EES}),
let us prove the class number formula (\ref{ECNF}). 

\begin{proof}[Proof of Theorem \ref{T2}]
For primes $l \ne p$, equation (\ref{ECNF}) claims that the $l$-class 
number of $L$ is given by $h_l(L) = (h_l(K)/h_l(F))^2 h_l(k)$; in fact 
we have an isomorphism
\begin{equation}\label{ECI}
  \Cl_l(L) \simeq \Cl_l(k) \times \Cl_l(K/F) \times \Cl_l(K'/F),
\end{equation}
where $\Cl(K/F)$ is the relative class group of $K/F$ defined
by the exact sequence
$$ \begin{CD} 
   1 @>>> \Cl(K/F) @>>> \Cl(K) @>{N_{K/F}}>>  \Cl(F) @>>> 1;
   \end{CD} $$
note that the norm is surjective since $K/F$ is nonabelian of
prime degree.

The isomorphism (\ref{ECI}) follows from the fact that the transfer 
of ideal classes $j_{k \to L}:$ is injective and the norm $N_{L/k}$ 
is surjective on classes of order coprime to $p$, hence 
$N_{L/k} \circ j_{k \to L}(c) = c^l$ induces an automorphism of 
$\Cl_l(k)$, which in turn implies that the sequence
$$ \begin{CD} 
   1 @>>> \Cl_l(L)[N] @>>> \Cl_l(L) @>>> \Cl_l(k) @>>> 1 
   \end{CD} $$
splits, i.e. $\Cl_l(L) \simeq \Cl_l(k) \times \Cl_l(L)[N]$.
A result of Jaulent (see \cite{CM}) guarantees that the
transfer of ideal classes $\Cl_2(K) \lra \Cl_2(L)$ is 
injective; for primes $l \nmid 2p$, the injectivity of 
$\Cl_l(K) \lra \Cl_l(L)$ is trivial. By Furtw\"angler's 
Theorem 90 and Lemma \ref{LST} we have
$\Cl_l(L)[N] = \Cl_l(K/F) \Cl_l(K'/F)$. 

We claim that $\Cl_l(K/F) \cap \Cl_l(K'/F) = 1$: a class
$c \in \Cl_l(K/F) \cap \Cl_l(K'/F)$ is fixed by $\tau$ and
$\sigma\tau$, hence by $\sigma$, and since it is killed by 
the norm, we find $c^p = 1$; since $c$ has $l$-power order,
this implies $c = 1$.

It remains to prove the $p$-part of the class number formula.
In the rest of the proof, all class numbers and class groups
are $p$-class numbers and $p$-class groups.

Let $N =N_{L/k}$ denote the relative norm of $L/k$. Since $L/k$
is unramified and cyclic, we know that $(\Cl(k):N\Cl(L)) = p$.
Thus 
\begin{equation}\label{ECNF1}
h_L = \# \Cl(L) = \# \Cl(L)[N] \cdot \# N\Cl(L) 
 = \frac{h_K^2}{\# \Am^{1+\tau}} \cdot \frac{h_k}p
 = p^{e-\rho} h_K^2 h_k. \end{equation}

If $B$ is a subgroup of finite index in an abelian group $A$
and if $f: A \lra A'$ is a group homomorphism, then
$(A:B) = (A^f:B^f)(\ker f: \ker f \cap B)$, where $A^f$ and
$B^f$ denote the images of $A$ and $B$ under $f$.

Now let us apply this to the special situation where $f$ is
given by the norm map $N: E_L/E_KE_{K'}E_k \lra E_k/E_k^pNE_K$. We have
$$ (E_L:E_KE_{K'}E_k) = (NE_L:N(E_KE_{K'}E_k)) \cdot
                        (E_L[N]: E_L[N] \cap E_KE_{K'}E_k).$$

\begin{lem}
If $L/k$ is unramified, then $(E_L[N]: E_L[N] \cap E_KE_{K'}E_k) = 1$.
\end{lem}

\begin{proof}
It suffices to show that $E_L[N] \subseteq E_K E_{K'}$. Assume 
therefore that $N_{L/k}\eps = 1$ for some $\eps \in E_L$. Then 
$\eps = \alpha^{1-\sigma}$ for some $\alpha \in L^\times$ by
Hilbert's Theorem 90, hence $\fa = (\alpha)$ is ambiguous. Since
$L/k$ is unramified, $\fa$ must be an ideal from $k$, and this 
implies that $\alpha = \eta a$ for some $\eta \in E_L$ and 
$a \in k^\times$. But then $\eps = \alpha^{1-\sigma} = \eta^{1-\sigma}
\in E_L^{1-\sigma}$, and by Lemma \ref{LST} we have
$E_L^{1-\sigma} = (E_L^{-\sigma})^{1+\tau} E_L^{1+\sigma\tau}
  \subseteq E_K E_{K'}$.
\end{proof}

Thus $q = (NE_L:N(E_KE_{K'}E_k))$; clearly 
$N(E_KE_{K'}E_k) = E_k^pNE_K$, and we can transform $q$ as follows:
\begin{align*}
(NE_L:N(E_KE_{K'}E_k)) 
          & =   (NE_L:E_k^pNE_K) 
          \ = \ \frac{(E_k:E_k^pNE_K)}{(E_k:NE_L)} \\ 
          & = \frac{(E_k:E_k^p)(E_k^p:E_k^pNE_K)}{(E_k:NE_L)} \\
          & = \frac{(E_k:E_k^p)(E_k^pE_F:E_k^pNE_K)}
                   {(E_k^pE_F:E_k^p)(E_k:NE_L)}.
\end{align*}
Now
\begin{align*}
(E_k^pE_F:E_k^pNE_K) & = \frac{(E_k^pE_F : E_k^p)}{(E_k^pNE_K : E_k^p)} 
       \ = \ \frac{(E_F:E_F \cap E_k^p)}{(NE_K : NE_K \cap E_k^p)} \\
       & = \frac{(E_F:E_F^p)}{(NE_K:E_F^p)} \ = \ (E_F:NE_K), 
\end{align*}
as well as
$$ (E_k^pE_F:E_k^p) = (E_F : E_F \cap E_k^p) = (E_F : E_F^p),$$
hence we get
$$ q \ = \ 
   \frac{(E_k:E_k^p)(E_F:NE_K)}{(E_F : E_F^p)(E_k:NE_L)} 
   \ = \ p^{\lambda(k) - \lambda(F)} p^{e+1-\rho}, $$
where $\lambda(M)$ denotes the $\Z$-rank of the unit group of a 
number field $M$. Note that $W_L = W_k$ (where $W_M$ denotes the 
group of roots of unity in $M$) since $L/F$ is non-abelian. 

Collecting everything we find

$$ h_L  = p^{e-\rho} h_K^2 h_k  = p^{-a} q h_K^2 h_k $$
for the $p$-class numbers, and this proves the theorem.
\end{proof}

\section{The Lower Bound for $r_p(K)$}

The idea of the proof is to lift parts of $\Cl_p(k)$ and 
$\Cl_p(K)$ to $L$ and compare their images. We start with the
group $\Cl(k)[p]$ of rank $r_p(k)$; its image after lifting it 
to $\Cl(L)$ has rank $\rk \Cl(k)[p]^j = r_p(k)-\rho$. Now 
observe that $\Cl(k)[p]^j$ is a subgroup of $\Cl_p(L)$ that is 
killed by $p$, $1+\tau$, $\sigma-1$, and the relative norm 
$N = N_{L/k}$. In particular, $\Cl(k)[p]^j \subseteq \cC_0$, where 
$\cC_0  = \{c \in \Cl_p(L): c^p = c^{1+\tau} = c^{1-\sigma} = 1\}$.
The key result is the following observation:

\begin{prop}\label{Pk}
There exists a monomorphism $\cC_0 \hra \Cl(K)[p]/\Am^+$.
\end{prop}

We know that 
$\rk \Cl(K)[p] = r_p(K)$ and $\rk \Am^+ = \rho-1-e$; 
since both groups are elementary abelian we deduce that 
$\rk \Cl(K)[p]/\Am^+ = r_p(K) - \rho + e + 1$. Thus 
from $\Cl(k)[p]^j \subseteq \cC_0 \subseteq \Cl(K)[p]/\Am^+$
we deduce that  
$$ r_p(k) - \rho = \rk \Cl(k)[p]^j \le \rk \Cl(K)[p]/\Am^+ 
                 =  r_p(K) - \rho + e + 1, $$ 
and this proves Theorem \ref{T1}.

It remains to prove Proposition \ref{Pk}. The next result 
(showing in particular that $\Am^+ \subseteq \Cl_p(K)$) 
can be found in Halter-Koch \cite{HK}:

\begin{lem}
Let $L/F$ be as above; in particular, assume that $L/k$ is 
unramified. We have $\Cl_p(L)[N] = \Cl_p(K)\Cl_p(K')$ and
$\Cl_p(K) \cap \Cl_p(K') = \Am^+$, where $K$ and $K'$ are
the fixed fields of $\tau$ and $\sigma\tau$.
\end{lem}

\begin{proof}
Since $(L:K)= 2$, the transfer of ideal classes $\Cl_p(K) \lra \Cl_p(L)$
is injective, and we can view $\Cl_p(K)$ as a subgroup of $\Cl_p(L)$.
Clearly $\Cl_p(K)$ and $\Cl_p(K')$ are killed by $N$, so 
$\Cl_p(K)\Cl_p(K') \subseteq \Cl_p(L)[N]$.

Using Lemma \ref{LST} we now find
$$\Cl_p(L)^{1-\sigma} = \Cl_p(L)^{1+\tau}\Cl_p(L)^{1+\sigma\tau}
  \subseteq  \Cl_p(K)\Cl_p(K') \subseteq \Cl_p(L)[N], $$ 
and by Furt\-w\"ang\-ler's Theorem 90 we have equality throughout.
\end{proof}

\begin{proof}[Proof of Prop. \ref{Pk}]
Given any $c \in \Cl_p(L)[N]$, we can write $c = c_1c_2$
with $c_1 \in \Cl_p(K)$ and $c_2 \in \Cl_p(K')$. Since
$\Cl_p(K) \cap \Cl_p(K') = \Am^+$, the $c_i$ are determined
modulo $\Am^+$, and we get a homomorphism 
$\lambda: \Cl_p(L)[N] \lra \Cl_p(K)/\Am^+$. 

We claim that $c_1 \in \Cl(K)[p]$ if $c \in \cC_0$. To prove
this, assume that $c^\sigma = c$ and $c^p = 1$; from 
$c_2^{\sigma\tau} = c_2$ we get $c_2^\sigma = c_2^\tau$,
and since $c^\tau = c^{-1}$ and $c_1^\tau = c_1$, we find
$c^{-1} = c^\tau = c_1 c_2^\tau$, that is, $c_2^\tau = c_1^{-2}c_2^{-1}$.
This gives $c_1^\sigma = (cc_2^{-1})^\sigma = c_1c_2c_1^2c_2 = c_1^3 c_2^2$.
Induction shows that $c_1^{\sigma^t} = c_1^{2t+1} c_2^{2t}$ and
$c_2^{\sigma^t} = c_1^{-2t} c_2^{1-2t}$. In particular,
$$ c_1^\nu \ = \ c_1^{1 + 3 + 5 + \ldots + 2p-1} c_2^{2 + 4 + \ldots + 2p-2} 
    \ = \  c_1^{p^2} c_2^{(p-1)p} = c^{(p-1)p} c_1^p.$$
But since $c^p = 1$ and $c_1^\nu = 1$, this implies $c_1^p = 1$,
that is, $c_1 \in \Cl(K)[p]$.

Now assume that $c \in \ker \lambda$; then $c_1 \in \Am^+$,
hence $c = c_1c_2 \in \Cl(K')$. Thus $c$ is fixed by $\sigma$ and
$\sigma\tau$, hence by $\tau$; since $c \in \Cl(L)^-$, this 
implies $c^2 = c^{1+\tau} = 1$, hence $c = 1$. Thus $\lambda$
is injective.
\end{proof}

\section{Embedding Problems}
Theorem \ref{TNo} on the existence of unramified $E(p^3)$-extensions
is a rather simple consequence of our results. Let $k/F$ be a
quadratic extensions, $p$ and odd prime such that $p \nmid h(F)$,
and $L/F$ a normal extension with $\Gal(L/F) \simeq D_p$ and $L/k$
unramified. If $\rk \Cl_p(k) \ge 2+e$, then Theorem \ref{T1} 
guarantees that any nonnormal subextension $K$ of $L/F$ will have 
class number divisible by $p$. Let $M/K$ be an unramified cyclic
extension of $K$, and let $N$ denote the normal closure of $LM/k$.
Then $N/k$ is a $p$-extension containing $L$, and its maximal 
abelian subextension $N^{\text{ab}}$ has type $(p,p)$. Let $E/k$ 
be a central extension of $N^{\text{ab}}/k$ of degree $p^3$ over $k$; 
then $\Gal(E/k) = E(p^3)$ or $\Gal(E/k) = \Gamma(p^3)$, where 
$\Gamma = \Gamma(p^3)$ is the nonabelian group of order $p^3$ 
and exponent $p^2$. We claim that $\Gal(E/k) = E(p^3)$. 

We remark in passing that -- in the case where $N \ne E$ --  
this follows immediately from the fact that $\Gamma(p^3)$ has 
trivial Schur multiplier. In general, we have to invoke the 
automorphism group $\Aut(\Gamma)$ of $\Gamma(p^3)$. It is 
known (see Eick \cite{Eick} and e.g. Schulte \cite{Schul}) that 
$\Aut(\Gamma) \simeq (\Z/p\Z)^\times \times p\text{-group}$,
and that $(\Z/p\Z)^\times$ acts trivially on exactly one of 
the two generators of $\Gamma/\Gamma'$. In particular, the 
unique element of order $2$ in $\Aut(\Gamma)$ acts as $-1$ 
on one and trivially on the other generator.

On the other hand, $\Gal(N^{\text{ab}}/F)$ is a generalized 
dihedral group by class field theory (since $p$ does not divide
the class number of $F$), hence the element of order $2$ in
$\Gal(k/F)$ acts as $-1$ on both generators of $\Gal(N^{\text{ab}}/k)$:
this means that $\Gal(N/k) \ne \Gamma(p^3)$, and Theorem \ref{TNo} follows.

\section{The Upper Bound for $r_p(K)$}

We will start by proving the upper bound in Theorem \ref{T1}
in two special cases: a refinement of our techniques used to derive
the lower bound will give the result if $p = 3$, and a simple Galois
theoretic argument suffices to prove it in the case $r_p(k) = 1$.

\subsection{The Case $p = 3$} \hfill \\

\noindent
In the special case $p = 3$ we can prove the upper bound
$r_p(K) \le r_p(k)-1$ using the same techniques we used
for deriving the lower bound. Our first ingredient holds in
general:

\begin{lem}
We have $\rk \Am_\st[N] = r_p(k) - \rho$ and 
$\rk \Am^-[N] \le r_p(k) - \rho + e$.
\end{lem}

\begin{proof}
Clearly $\Cl(k)[p]^j \subseteq \Am_\st[N]$; conversely, if
$c \in \Cl_p(k)$ with $c^j \in \Am_\st[N]$, then $c^p = 1$,
hence we actually have $\Cl(k)[p]^j = \Am_\st[N]$, and this 
proves that $\rk \Am_\st[N] = r_p(k) - \rho$.

The exact sequence 
$$ \begin{CD}
   1 @>>> \Am_\st[N] @>>> \Am^-[N] @>>> E_F/NE_K \end{CD} $$
derived from (\ref{EES}) shows that  $\rk \Am^-[N] \le \rk \Am_\st[N] + e$.
\end{proof}

From now on assume that $p = 3$; then the map $c \too c^{1+2\sigma}$ 
defines a homomorphism $\mu: \Cl(K)[p] \lra \cC_0$.
In fact, since $c^\nu = c^3 = 1$, we have 
$\mu(c)^\sigma = c^{\sigma+2\sigma^2} = c^{-2-\sigma} = \mu(c)$
since $c^{\sigma^2} = c^{-1-\sigma}$. Moreover, 
$\mu(c)^\tau = c^{\tau(1 + 2\sigma^2)} = c^{-1-2\sigma} = c^{-1}$ 
since $c^\tau = c$ for $c \in \Cl(K)$; thus $\mu(c)$ is killed by $N$, 
$p$ and $1+\tau$, hence $\mu(c) \in \cC_0$. 

Next $c \in \ker \mu$ implies $c = c^\sigma$, i.e., $c \in \Am^+$,
and clearly $\Am^+ \subseteq \ker \mu$: thus

\begin{prop}
If $p = 3$, then $\Cl(K)[p]/\Am^+ \simeq \cC_0$.
\end{prop}

In particular, $r_p(K) - \rho + 1 + e = \rk \cC_0$ if $p = 3$. Since 
$\cC_0 \subseteq \Am^-[N]$, we find  
$r_p(K) - \rho + 1 + e \le r_p(k) - \rho + e$, and we have proved

\begin{thm}
If $p = 3$, then $r_p(k)-1-e \le r_p(K) \le r_p(k)-1$.
\end{thm}

\subsection{The case $r_p(k) = 1$} \hfill \\

\noindent
The second special case of the upper bound that can be proved easily is

\begin{prop}
If $r_p(k) = 1$, then $r_p(K) = 0$.
\end{prop}

\begin{proof}
Assume not; then there exists a cyclic unramified extension $M/K$ 
of degree $p$. Let $N$ denote the normal closure of $ML/k$. If 
$N = ML$, then $ML/k$ has a Galois group of order $p^2$ and thus
is abelian, and since $ML/L$ and $L/k$ are unramified, so is $ML/k$.
Since $\Cl_p(k)$ is cyclic by assumption, we conclude that 
$\Gal(ML/k) = \Z/p^2\Z$, and since $p$ does not divide the class
number of $F$, we conclude that $ML/F$ is normal and
$\Gal(ML/F) \simeq D_{p^2}$. On the other hand, 
$\Gal(ML/K) \simeq \Z/2\Z \times \Z/p\Z$ by construction, and
since the dihedral group of order $2p^2$ does not contain an 
abelian subgroup of order $2p$, we have a contradiction.  
\end{proof}

\subsection{The Case $\Cl_p(k) \simeq (p,p)$}

Our main tool will be the following result
due to G.~Gras \cite{Gras}:

\begin{prop}\label{PGr}
Let $p$ be an odd prime, $G = \la \sigma \ra$ a group 
of order $p$ generated by $\sigma$, and assume that $G$ 
acts on the abelian $p$-group $A$ in such a way that 
$\# A^G = \# \{a \in A: a^{\sigma-1} = 1\} = p$. Put 
$\nu = 1 + \sigma + \sigma^2 + \ldots + \sigma^{p-1}$, 
let $n$ be the smallest positive integer such that
$A^{(\sigma-1)^n} = 1$, and write $ n = \alpha(p-1) + \beta$ 
with $0 \le \beta \le p-2$. 

If $A^\nu = 1$, then
$$ A \simeq (\Z/p^{\alpha+1}\Z)^\beta \times (\Z/p^\alpha\Z)^{p-1-\beta}. $$ 

If $A^\nu \ne 1$,  then 
$$ A \simeq \begin{cases}
     (\Z/p^2\Z) \times (\Z/p\Z)^{n-2} & \text{if}\ n < p, \\
     (\Z/p\Z)^p                       & \text{if}\ n = p, \\ 
     (\Z/p^{\alpha+1}\Z)^\beta \times (\Z/p^\alpha\Z)^{p-1-\beta} 
                                      & \text{if}\ n > p. 
    \end{cases} $$
Note that $\# A = p^n$ and that the $p$-rank of $A$ is bounded
by $p$.
\end{prop}

Assume now that $\Cl_p(k) \simeq (p,p) = \Z/p\Z \times \Z/p\Z$
and put $A = \Cl_p(L)$. Then $\# A^G = \# \Am(L/k) = p$ by the 
ambiguous class number formula since every unit in $k$ is a norm 
from $L$. Moreover, $A_k = N_{L/k} A$ has index $p$ in $\Cl_p(k)$ by class 
field theory, hence $A_k = \la c \ra$ for some ideal class $c$
of order $p$, and we have to distinguish two cases:
\begin{enumerate}
\item[(A)] $c$ capitulates in $L/k$; then $A^\nu = 1$.
\item[(B)] $c$ does not capitulate in $L/k$; then $A^\nu \ne 1$.
\end{enumerate}
Moreover, $\rho \ge e+1$ implies that the following classification is 
complete:
$$ \# \Cl_p(K) \cap \Cl_p(K') = \begin{cases}
      1 & \text{if}\ (\rho, e) = (1, 0), (2,1) \\
      p & \text{if}\ (\rho, e) = (2,0). \end{cases} $$
Applying the class number formula (\ref{ECNF1}) we get
$$ h_p(L) = p^{2+e-\rho} h_K^2 = p^\mu $$
for some integer $\mu \ge 1$. Now we can prove

\begin{thm}
Let $p$ be a prime and assume that $F$ is a number field whose 
class number is not divisible by $p$. Let $L/F$ be a normal 
extension with Galois group $D_p$, and let $L/k$ be unramified. 
Assume that $\Cl_p(k) \simeq (\Z/p\Z)^2$. Then
\begin{enumerate}
\item[a)] $h_p(L) = h_p(K)^2 p^{2-\rho} = p^\mu$ for some $\mu \ge 2$,
      and in particular $\mu \equiv \rho \bmod 2$. 
\item[b)] Write $\mu = \alpha(p-1) + \beta$ with $0 \le \beta < p-1$;
      then the structure of $\Cl_p(L)$ is given by the following table:
      $$ \begin{array}{|c||c|c|} \hline
        \text{case} & A & B \\ \hline\hline
  \rsp  \mu > p &  \multicolumn{2}{c|}{(\Z/p^{\alpha+1}\Z)^\beta \times
                     (\Z/p^\alpha\Z)^{p-1-\beta}} \\ \hline
  \rsp  \mu = p & (\Z/p^2\Z) \times (\Z/p\Z)^{p-2} & (\Z/p\Z)^p \\ \hline
  \rsp  \mu < p & (\Z/p\Z)^\mu & (\Z/p^2\Z) \times (\Z/p\Z)^{\mu-2} \\ \hline 
       \end{array} $$
\item[c)] Write $\mu-1 = a(p-1) + b$, $0 \le b \le p-2$; note that 
      $b$ is even if $\rho-e$ is odd. The structure of $\Cl_p(L)[N]$
      and  $\Cl_p(K)$ is given by
      \begin{align*} \quad 
       \Cl_p(L)[N] & \simeq (\Z/p^{a+1}\Z)^b \times (\Z/p^a\Z)^{p-1-b}, \\
       \Cl_p(K) & \simeq \begin{cases}
         (\Z/p^{a+1}\Z)^{b/2} \times  (\Z/p^a\Z)^{(p-1-b)/2}
                & \text{if}\ \rho = e+1, \\
         (\Z/p^{a+1}\Z)^{(b+1)/2} \times  (\Z/p^a\Z)^{(p-2-b)/2}
                & \text{if}\ \rho = 2, e = 0. \end{cases} 
       \end{align*}
      Observe that $\Cl_p(K)$ is elementary abelian if and only if 
      $\mu \le p$. On the other hand, we have $\rk \Cl_p(K) = \frac{p-1}2$ 
      whenever $\mu \ge p-1$.
\end{enumerate}
\end{thm}

\begin{proof}
We already proved the class number formula in a), and b) follows
by applying Prop. \ref{PGr} to $A = \Cl_p(L)$. Similarly, the claims 
in c) about the structure of $\Cl_p(L)[N]$ follow by applying
Prop. \ref{PGr} to $A = \Cl_p(L)[N]$. 

It remains to derive the structure of $\Cl_p(K)$. If $\rho = e+1$,
then we have seen that $\Cl_p(K) \cap \Cl_p(K') = 1$, hence
$\Cl_p(L)[N] \simeq \Cl_p(K) \oplus  \Cl_p(K)$, and this allows us
to deduce the structure of $\Cl_p(K)$ from that of $\Cl_p(L)$.
If $(\rho,e) = (2,0)$, on the other hand, then 
$\Cl_p(L)[N] = \Cl_p(K) \Cl_p(K')$ with 
$\Cl_p(K) \cap \Cl_p(K') \simeq \Z/p$, and again the claims follow
easily.
\end{proof}

\medskip\noindent{\bf Examples.}
Consider the cubic field $K_a$ generated by a root of the
polynomial $x^3 + ax + 1$; let $d = \disc k$. 
$$ \begin{array}{r|r|r|c|c}
  a & d \ \ & \Cl_3(k) & \Cl_3(K) & \Cl_3(L) \\ \hline
 29 &      -97583 & (3,3) &   (3) & (3,3,3) \\
 10 &       -4027 & (3,3) &   (3) & (3^2,3) \\
 70 &    -1372027 & (3,3) & (3^2) & (3^3,3^2) \\
 94 &    -3322363 & (3,3) & (3^3) & (3^4,3^3) \\
755 & -1721475527 & (3,3) & (3^4) & (3^5,3^4) \\
409 &  -273671743 & (3,3) & (3^5) & (3^6,3^5) 
 \end{array} $$ 
The data suggest that the exponent of $\Cl_3(K_a)$ is not bounded.


\section{Examples}

It is expected that the upper bounds are best possible
even if $p > 3$. The following family of simplest dihedral 
quintics extracted from Kondo \cite{Kon} show that the upper 
bound is attained for $p=5$. Let
$$ f(x) = x^5-2x^4+(b+2)x^3-(2b+1)x^2+bx+1, $$
let $\alpha$ denote a root of $f$, and put $K = \Q(\alpha)$.
Then $\disc K = d^2$ for some odd $d$, and if we choose
the sign of $d$ such that $d \equiv 1 \bmod 4$, then
then the splitting field $L$ of $f$ (which has Galois group $D_5$)
is unramified over its quadratic subfield $k = \Q(\sqrt{d}\,)$ if 
$d$ is squarefree. 
$$ \begin{array}{r|r|c|c|c}
 b & d \ \  & \Cl(k) & \Cl(K) & \Cl(L) \\ \hline
 1 &    -103 & (5)     & 1  & 1 \\
19 &  -38047 & (15,5)  & (20,4) & (300,20,4,4) \\
39 & -280847 & (20,20) & (55,5) & (1100, 220, 5, 5) 
 \end{array} $$

Using $F = \Q(\sqrt{5}\,)$ as the base field, we find
$$ \begin{array}{r|r|r|r}
 b & d & \Cl(k) & \Cl(K) \\ \hline
41 &     47 & (5)         & 1  \\
 9 &   5447 & (60, 20)    & (55) \\
16 &  23983 & (50, 10, 5) & (305) \\
17 &  28199 & (480, 15)   & (4, 4) \\
39 & 280847 & (1080,40)   & (55,5) 
 \end{array} $$ 
 
Here are a few examples for $p = 3$ that also show that the term 
$e$ in our lower bound is necessary: let $d$ be the discriminant of 
a dihedral cubic number field $k_0$, and consider the fields
$F = \Q(\sqrt{-3}\,)$, $K = k_0(\sqrt{-3}\,)$,
$k = \Q(\sqrt{-3},\sqrt{d}\,)$ and $L = Kk$. 

$$ \begin{array}{r|r|r}
   d   &  \Cl(k) & \Cl(K) \\ \hline
      -31 &  (3)      & (1) \\
     -107 &  (3,3)    & (1) \\
    -4027 &  (3,3,3)  & (6,2) \\
    -8751 &  (12,3,3) & (3,3) \\   
      229 &  (6,3)    & (2) \\
      469 &  (6, 6)   & (3) \\
    26821 &  (72,3)   & (18) \\ 
  2813221 &  (198,6,6,6) & (285,3) \\ 
 13814533 & (270,3,3,3,3) & (360,3,3) 
   \end{array} $$

\end{document}